\documentclass[10pt, a4paper]{article}
\usepackage{amscd}
\usepackage{amsfonts}
\begin{document}

\newtheorem{theorem}{Theorem}
\newtheorem{corollary}[theorem]{Corollary}
\newtheorem{definition}[theorem]{Definition}
\newtheorem{lemma}[theorem]{Lemma}
\newtheorem{proposition}[theorem]{Proposition}
\newtheorem{remark}[theorem]{Remark}
\newtheorem{example}[theorem]{Example}
\newtheorem{notation}[theorem]{Notation}
% fim da demonstracao
\def\Qed{\hfill\raisebox{.6ex}{\framebox[2.5mm]{}}\\[.15in]}

\title{On equations of double planes with $p_g=q=1$}

\author{Carlos Rito}

\date{}
\pagestyle{myheadings}
\maketitle
\setcounter{page}{1}

\begin{abstract}

This paper describes how to compute equations of
plane models of minimal Du Val double planes of
general type with $p_g=q=1$ and $K^2=2,\ldots,8.$
A double plane with $K^2=8$ having bicanonical map not composed
with the associated involution is also constructed.
The computations are done using the algebra system Magma.

\noindent 2000 Mathematics Classification: 14J29, 14Q05.
\end{abstract}

\section{Introduction}

There are two typical methods to construct examples of surfaces:
Campedelli --- double covers ramified over curves, usually highly singular
--- and Godeaux --- quotients by group actions ({\it cf.} \cite{Re}).
In this paper the first method is used to obtain
new examples of double planes of general type with $p_g=q=1.$
Often the singularities impose too many conditions
on the linear parameters of systems of curves, and this implies hard calculations.
The {\em Computational Algebra System MAGMA} (\cite{BCP})
is used to perform the computations.

Several authors have studied surfaces of general type with $p_g=q=1$
(\cite{Ca1}, \cite{Ca2}, \cite{CC1}, \cite{CC2}, \cite{CP},
\cite{Po1}, \cite{Po2}, \cite{Po3}, \cite{Po4}, \cite{Pi}),
but these surfaces are still not completely understood, and few examples are known.

For such a smooth minimal surface $S,$ one has $2\leq K^2\leq 9$ and $K^2=9$ only if
the bicanonical map of $S$ (given by $|2K|$) is birational (from \cite{CM} and \cite[Th\' eor\` eme 2.2]{Xi1}).
Since $q=1,$ the Albanese map of $S$ is a connected fibration onto an elliptic curve.
We denote by $g$ the genus of a general Albanese fibre of $S.$

Du Val (\cite{Du}) classified the regular surfaces $S$ of general type
with $p_g\geq 3$ whose general canonical curve is smooth
and hyperelliptic.  For these surfaces the bicanonical map is
composed with an involution $i$ such that $S/i$ is rational.
The families of surfaces exhibited by Du Val
are nowadays called the {\em Du Val examples}.

The {\em standard case} of non-birationality of the bicanonical morphism $\phi_2$
is the case where $S$ has a genus 2 fibration.

For the non-standard case with $\phi_2$ of degree 2 and bicanonical
image a ruled surface, Xiao Gang (\cite[Theorem 2]{Xi2}) extended Du Val's list to
$p_g(S)\geq 1$ and added two extra families (this result is still valid
assuming only that $\phi_2$ is composed with an involution such that
the quotient surface is a ruled surface).
Later G. Borrelli (\cite{Br}) excluded these two families,
confirming that the only possibilities for this case are the Du Val examples.

This way, if the bicanonical map $\phi_2$ of $S$ is composed with $i,$ $S/i$ is ruled
and $S$ presents the non-standard case, then $S$ is a Du Val double plane (see Definition \ref{DV}).
But, to my knowledge, the existence of these double planes with $p_g=q=1$ and
$K^2=4,\ldots,7$ has not been shown yet, and no equation has been given for $K^2=3$ or $8.$

In this paper we compute equations of plane models of minimal Du Val double planes of general type
with $p_g=q=1$ and $K^2=2,\ldots,8.$ We also construct a double plane with $K^2=8$ having
bicanonical map not composed with the associated involution.

Notice that we give the first example of a minimal surface of general type with $p_g=q=1$ and $K^2=7.$
We also give the first example with $K^2=5$ having Albanese fibration of genus $g\ne 2.$

The first example with $K^2=6$ has been given by the author in \cite{Ri}, as a double cover of a $K3$ surface.
Polizzi (\cite{Po2}) also gives examples with $K^2=6,4$ and $g\ne 2.$
These surfaces contain $8-K^2$ $(-2)$-curves, all in the same Albanese fibre.
One can verify that, in the double plane examples with $K^2=4$ or $6,$
there is no Albanese fibre containing $8-K^2$ $(-2)$-curves.
Therefore our surfaces are different from the above examples.

In \cite{Po4}, Polizzi classifies surfaces of general type with $p_g=q=1,$ $K^2=8$ and
bicanonical map of degree 2. He constructs examples using quotients under the action of a group
and shows that these surfaces are Du Val double planes, describing the branch locus of
the respective plane models.
In sections \ref{K8642} and \ref{K8g4} we show how to obtain equations for such branch loci.

Section \ref{prelim} contains a structure description of the Albanese fibration of $S$
and the definition of Du Val double plane.
In Section \ref{ExplDblPlns} we describe the principal steps of the constructions.
To obtain the equations of the branch loci we impose conditions to the linear
parameters of systems of plane curves. Since the ramification curves are contained
in Albanese fibres, in some cases it is easier to start by
constructing non-reduced Albanese fibres, which simplify the computations.
These are done using the algebra system Magma in Appendix \ref{Magma}.

\bigskip
\noindent{\bf Notation and conventions}

We work over the complex numbers;
all varieties are assumed to be projective algebraic.
For a projective smooth surface $S,$ the {\em canonical class} is denoted by
$K,$ the {\em geometric genus} by $p_g:=h^0(S,\mathcal O_S(K)),$ the {\em irregularity}
by $q:=h^1(S,\mathcal O_S(K))$ and the {\em Euler characteristic} by
$\chi=\chi(\mathcal O_S)=1+p_g-q.$

An {\em $(-n)$-curve} $C$ on a surface is a curve
isomorphic to $\mathbb P^1$ such that $C^2=-n$. We say that a curve
singularity is {\em negligible} if it is either a double point or a
triple point which resolves to at most a double point after one
blow-up.
An $(m_1,m_2,\ldots)$-point, or point of order $(m_1,m_2,\ldots),$
is a point of multiplicity $m_1,$ which resolves to a point
of multiplicity $m_2$ after one blow-up, etc.

An {\em involution} of a surface $S$ is an
automorphism of $S$ of order 2. We say that a map is {\em composed with
an involution} $i$ of $S$ if it factors through the double cover $S\rightarrow
S/i.$

The rest of the notation is standard in algebraic geometry.

\bigskip
\noindent{\bf Acknowledgements}

The author is a collaborator of the Center for Mathematical
Analysis, Geometry and Dynamical Systems of Instituto Superior T\'
ecnico, Universidade T\' ecnica de Lisboa, and
is a member of the Mathematics Department of the
Universidade de Tr\'as-os-Montes e Alto Douro.
This research was partially supported by FCT (Portugal) through
Project POCTI/MAT/44068/2002.

\section{Preliminaries}\label{prelim}

We say that a smooth surface $S$ is a {\em double plane} if $S$ has an involution $i$ such
that $S/i$ is a rational surface. A {\em plane model} of $S$ is a double cover $X\rightarrow\mathbb P^2$
such that $X$ is a normal surface and there exists a commutative diagram
$$
\begin{CD}\ S@> >>X\\ @V VV  @VV  V\\ S/i@> >> \mathbb P^2
\end{CD}
$$
such that the horizontal arrows denote birational maps.

\subsection{Double planes with $q=1$}\label{fA}

Let $S$ be a smooth minimal surface of general type with an
involution $i.$ Since $S$ is minimal of general type, this
involution is biregular. The fixed locus of $i$ is the union of a
smooth curve $R''$ (possibly empty) and of $t\geq 0$ isolated points
$P_1,\ldots,P_t.$ Let $S/i$ be the quotient of $S$ by $i$ and
$p:S\rightarrow S/i$ be the projection onto the quotient. The
surface $S/i$ has nodes at the points $Q_i:=p(P_i),$ $i=1,\ldots,t,$
and is smooth elsewhere. If $R''\not=\emptyset,$ the image via $p$
of $R''$ is a smooth curve $B''$ not containing the singular points
$Q_i,$ $i=1,\ldots,t.$ Let now $h:V\rightarrow S$ be the blow-up of
$S$ at $P_1,\ldots,P_t$ and set $R'=h^*(R'').$ The involution $i$
induces a biregular involution $\widetilde{i}$ on $V,$ whose fixed
locus is $R:=R'+\sum_1^t h^{-1}(P_i),$ and the quotient
$W:=V/\widetilde{i}$ is smooth.

Suppose now that $q(S)=1$ and $S/i$ is a rational surface.
Since $q=1,$ the Albanese variety of $S$ is an elliptic curve $E$ and the Albanese map is a
connected fibration (see {\it e.g.} \cite{Be} or \cite{BPV}).
This fibration is preserved by $i$ and so we have a commutative diagram
\begin{equation}\label{CmtDgr}
\begin{CD}\ V@>h>>S @>>>E\\ @V\pi VV  @V p VV @VV V\\ W@>>> S/i @>>>\mathbb P^1
\end{CD}.
\end{equation}
Denote by
$$f_A:W\rightarrow\mathbb P^1$$ the fibration induced by the Albanese
fibration of $S.$
The double cover $E\rightarrow\mathbb P^1$ is
ramified over 4 points $p_j$ of $\mathbb P^1,$ thus the branch locus
$B$ of $\pi:V\rightarrow W$ is contained in 4 fibres $$F_A^j:=f_A^*(p_j),\
j=1,...,4,$$ of the fibration $f_A$.
By Zariski's Lemma (see {\it e.g.} \cite{BPV}) the irreducible
components $B_i$ of $B$ satisfy $B_i^2\leq 0.$ 
Since  $\pi^*(F_A^j)$ has even multiplicity, each component of
$F_A^j$ which is not a component of the branch locus must be of even multiplicity.

In the following sections $W$ can always be contracted to $\mathbb P^2.$ We keep the same notation
for the image of $f_A$ and $F_A^j$ on $\mathbb P^2.$

\subsection{Du Val double planes}\label{defDuVal}

\begin{notation}\label{notation0}
Given distinct points $p_0,\ldots,p_j,\ldots,p_{j+s}\in\mathbb P^2,$ let $T_i$ denote the
line through $p_0$ and $p_i,$ $i=1,\ldots,j.$
We say that a plane curve is of type $$d\big(m,(n,n)_T^j,r^s\big)$$
if it is of degree $d$ and if it has: an $m$-uple point at $p_0,$
an $(n,n)$-point at $p_1,\ldots,p_j,$ an $r$-uple point at $p_{j+1},\ldots,p_{j+s}$
and no other non-negligible singularities.
The index $_T$ is used if $T_i$ is tangent to the $(n,n)$-point at $p_i.$

An obvious generalization is used if there are other singularities.
\end{notation}

\begin{definition}
Let $C_0$ and $F$ denote, respectively, the negative section and a rational fibre of
the Hirzebruch surface $\mathbb F_2.$

The {\em Bombieri-Du Val surface} is the minimal model of the double cover of $\mathbb F_2$
with branch locus a smooth curve in $C_0+|7C_0+14F|.$
\end{definition}

\begin{definition}\label{DV}
A {\em Du Val surface} is either
\begin{itemize}
  \item [$\mathcal{B}$)] the Bombieri-Du Val surface
\end{itemize}
or a minimal double plane having a plane model with branch locus $D$ one of the following:
\begin{itemize}
  \item [$\mathcal{D}$)] a smooth curve of degree 8;
  \item [$\mathcal{D}_0$)] a smooth curve of degree 10;
  \item [$\mathcal{D}_n$)] a curve of type $[10+2n]\big(2n+2,(5,5)_T^n\big),$ $n\in\{1,\ldots,6\};$
\end{itemize}
or one of $\mathcal{B},$ $\mathcal{D}$ or $\mathcal{D}_n,$
imposing additional $4$-uple or $(3,3)$-points to $D.$
\end{definition}

Surfaces of type $\mathcal{B},$ $\mathcal{D}$ or $\mathcal{D}_n$
are called {\em Du Val's ancestors}.
The Bombieri surface has $K^2=9,$ $p_g=6$ and $q=0;$
a Du Val ancestor of type $\mathcal{D}$ has $K^2=2,$ $p_g=3$ and $q=0;$
those of type $\mathcal{D}_n,$ $n\geq 0,$ have $K^2=8,$ $p_g=6-n$ and $q=0,$
except possibly in the case $n=6,$ where $p_g=q=1$ if the $6$ singular points
are contained in a conic.

Notice that the imposition of a 4-uple point to the branch locus decreases $K^2$ by 2 and
the Euler characteristic $\chi$ by 1, while a $(3,3)$-point decreases both $K^2$
and $\chi$ by 1. Negligible singularities in the branch locus do not change
these invariants.

The imposition of $6-n$ $4$-uple or $(3,3)$-points to the branch locus of a
Du Val surface of type $\mathcal{D}_n$ gives $p_g=0$ if the singular points
$p_1,\ldots,p_6$ are not contained in a conic and $p_g=q=1$ otherwise.

Notice that, in the case $\mathcal{D}_n,$ the branch locus contains the lines $T_1,\ldots,T_n.$
Each of these lines corresponds to two $(-2)$-curves in $W,$ which contract to two nodes of $S/i.$
A $(3,3)$-point in the branch locus also corresponds to a node of $S/i.$

The bicanonical map of a Du Val surface factors through a map of degree 2 onto a rational surface.
For more information on Du Val surfaces see \cite{Du} or \cite{Ci}.

\subsection{The Magma procedures $LinSys$ and $LinSys2$}\label{Proc}

Magma has the function $LinearSystem$ which calculates linear systems of plane curves with ordinary
singularities, but we want to work with non-ordinary singularities.
To achieve this, we define two procedures.
The first one, $LinSys,$ calculates the linear system $L$ of plane curves of
degree $d,$ in an affine plane $\mathbb A,$ having singular points $p_i$ of
order $(m1_i,m2_i)$ with tangent slope $td_i.$

The other procedure, $LinSys2,$ calculates the sub-system $J,$ of a given
linear system $L$ of plane curves, of those sections which have a singularity,
at a point $q,$ of type $m=(m_1,\ldots,m_j)$ with tangent
slopes given by $td=[td_1,\ldots,td_{j-1}].$

The code lines for these procedures are in Appendix \ref{Procedure}.

\section{Equations of plane models with $p_g=q=1$}\label{ExplDblPlns}

In this section we obtain equations of the branch locus $B\subset\mathbb P^2=\mathbb P(1,1,1)$ of
minimal Du Val double planes $S$ of type $\mathcal{D}_n$ (see Definition \ref{DV}) 
with $p_g=q=1$ and $K^2=2,\ldots,8.$ 
Notice that if $f=0$ is the defining equation of $B,$
then the equation $w^2=f$ gives a plane model of $S,$
in the weighted projective space $\mathbb P(({\rm deg\ }f)/2,1,1,1).$

In Section \ref{ExplBirrBic} we construct a double plane with $K^2=8$ having
bicanonical map not composed with the associated involution.

Table (\ref{tabular}) lists
the type of each branch curve that we are going to obtain and
the corresponding values of $(K^2,g),$
where $g$ denotes the genus of a general Albanese fibre of $S.$

We keep Notation \ref{notation0}.

\begin{equation}\label{tabular}
\begin{tabular}{c|c}

  Type of branch curve & $(K^2,g)$ \\ \hline
   & \\
  $22\big(14,(5,5)^6_T\big)$ & $(8,5),\ (8,4),\ (8,3)$ \\
   & \\
  $[10+2n]\big(2n+2,(5,5)^n_T,4^{6-n}\big)$ & $(6,4),\ (4,3),\ (2,2)$ \\
  $n=5,4,3$ &  \\
   & \\
  $[10+2n]\big(2n+2,(5,5)^n_T,(3,3),4^{5-n}\big)$ & $(7,5),\ (5,4),\ (3,3)$ \\
  $n=5,4,3$ & \\ 
   & \\
  $18\big(10,(5,5)^4_T,(3,3)^2\big)$ & $(6,3)$  \\
   & \\
  $26\big(14,(7,7)^5_T\big)$ & $(8,4)$  \\
\end{tabular}
\end{equation}
\\\\

First we find double planes with $(K^2,g)=(8,5)$ and $(K^2,g)=(7,5)$
such that the branch locus contains an element of the pencil $f_A$
which induces the Albanese fibration.
From here constructions with $(K^2,g)=(6,4),(5,4),(4,3),(3,3),(2,2)$ will follow easily.
A construction with $(K^2,g)=(6,3)$ and branch locus
strictly contained in elements of $f_A$ is also given.
Then we get surfaces with $(K^2,g)=(8,4)$ or $(8,3).$
These three surfaces with $K^2=8$ are surfaces of type I, II, III described by F. Polizzi in \cite{Po4}.
Finally we construct a pencil $l$ of plane curves such that the minimal double plane with branch locus
a general element of $l$ has $K^2=8$ and bicanonical map not composed with the corresponding involution.

In order to obtain $q=1,$ all the singularities, except the one at $p_0,$ are chosen
to be contained in a conic.

A difficulty that arises in the computation of singular curves is that often
the computer is not able to finish the calculations.
Except for Section \ref{K7}, our method is to try to figure out the configuration
of the non-reduced fibres $F_A^i$ of $f_A,$ described in Section \ref{fA}.
The problem of finding the support of $F_A^i$ deals
with curves of lower degree and with simpler singularities, hence
with faster computations. 

The next sections of this chapter describe the principal steps.
The detailed calculations are done in Appendix \ref{Magma},
using the Computational Algebra System Magma.

We keep Notation \ref{notation0}.
\subsection{$K^2=8,6,4,2$ and $g=5,4,3,2$}\label{K8642}
To construct a Du Val's ancestor $S$ of type $\mathcal{D}_6$ (it has $K^2=8$)
it suffices to find a curve $B'$ of type $$16\left(8,(4,4)^6_T\right),$$ singular at $p_0,\ldots,p_6,$
such that the branch locus $B:=B'+\sum_1^6T_i$ is reduced.
If $p_1,\ldots,p_6$ are contained in a conic, then $p_g(S)=q(S)=1.$

In this example the pencil $f_A$ which induces the Albanese fibration of $S$ has the
following elements:
$$F_A^1=2D_1+T_3+\cdots+T_6,\ \ \ F_A^2=2D_2+T_1+T_2,\ \ \ F_A^3=2D_3,\ \ \ F_A^4=D_4,$$
where $D_1,\ldots,D_4$ are curves of types
$$6\left(2,(2,2)^2_T,(2,1)^4_T\right),\ \ 7\left(3,(2,1)^2_T,(2,2)^4_T\right),\ \ 8\left(4,(2,2)^6_T\right),\ \ 16\left(8,(4,4)^6_T\right),$$
respectively.

To construct $B',$ we first find points $p_0,\ldots,p_6$ such that $D_1$ exists.
Then we use the procedure $LinSys$ (see Section \ref{Proc})
to verify the existence of $D_2$ and $D_3,$ also singular at $p_0,\ldots,p_6.$
The curve $B'=D_4$ is a general element of the pencil generated by $D_1, D_2, D_3.$\\

So let us look for $D_1.$ Briefly, the steps are as follows.
Let $\mathbb A$ be an affine plane, $C$ be a smooth conic not containing the origin $p_0$
of $\mathbb A$ and $p_1,\ldots,p_4$ be points in $C.$ Denote by $L$ the linear system of plane
curves of type 
$$6\left(2,(2,2)^2_T,(2,1)^2_T\right),$$
with singularities at $p_0,\ldots,p_4,$ respectively. Let $F$ be a general
element of $L$ and $p_5,p_6$ be general points of $\mathbb A.$ We define a scheme
$Sch$ by imposing the following conditions:
\begin{description}
  \item $\cdot$ $p_5,p_6\in C\bigcap F;$
  \item $\cdot$ $p_5,p_6$ are double points of $F;$
  \item $\cdot$ the singularities of $F$ at $p_5,p_6$ have one branch tangent to
  $T_5,T_6;$
  \item $\cdot$ $p_5\ne p_6$ and $p_5,p_6\not\in\{p_0,\ldots,p_4\}.$
\end{description}
Now we compute the points of $Sch$ with Magma, choosing one of the solutions
for $p_5,p_6,$ and we use the procedure $LinSys$ to compute $B',$ a general element
of the pencil of curves of type $16(8,(4,4)^6_T),$
with singularities at $p_0,\ldots,p_6.$

Finally we perform some verifications:
that ${B'}$ is reduced, the singularities are as expected, etc.\\

With this we find a minimal double plane with $p_g=q=1$ and $K^2=8.$
The divisor $2D_1+T_3+\cdots+T_i$ is also a good candidate for one
of the singular fibres ${F_A^i}$
in the case where the branch locus is a curve of type
$$[10+2n]\left(2n+2,(5,5)^n_T,4^{6-n}\right),\ \ \ n=5,4,3.$$
In fact, using the points $p_0,\ldots,p_6$ and the procedure $LinSys,$
one can find branch loci of those types, obtaining then minimal double
planes with $p_g=q=1,$ $K^2=6,4,2$ and $g=4,3,2,$ respectively.\\

The corresponding Magma calculations are in Appendix \ref{Appendix1}.
There we use symmetry in order to obtain faster computations.
\subsection{$K^2=7,5,3$ and $g=5,4,3$}\label{K7}
Here we impose a $(3,3)$-point to the branch locus of a Du Val's ancestor of type $\mathcal{D}_5,$
{\it i.e.} we construct a plane curve ${B'}$ of type $$15(7,(4,4)^5_T,(3,3))$$ such that
${B}:={B'}+\sum_1^5 T_i$ is reduced (see Notation \ref{notation0}).
The double cover with branch locus ${B}$ is a plane model
of a Du Val double plane $S$ with $K^2=7$ and $\chi=1.$

In this example the $(3,3)$-point $p_6$ is infinitely near to
the $(5,5)$-point $p_1$ of ${B}$
({\it i.e.} there is a $(5,5,3,3)$-point at $p_1$)
and the line $T_1,$ through $p_0,p_1,$
is tangent to the conic $C$ defined by $p_1,\ldots,p_5.$
This implies $q(S)=1,$ because $\widetilde{C}-\sum_1^6E_i$ is effective,
where the curves $E_i$ are the exceptional divisors corresponding to
the blow-ups at the points $p_i$ and $\widetilde C$ is the pullback of $C.$

In this case the pencil $f_A$ has elements
$$F_A^1=2D_1+T_4+T_5,\ \ \ F_A^2=2D_2+T_2+T_3,\ \ \ F_A^3=2D_3,\ \ \ F_A^4=D_4+T_1,$$
where $D_1,\ldots,D_4$ are curves of types
$$7\left(3,(2,2)^3_T,(2,1)^2_T,(2,2)\right),\ \ \ \ 7\left(3,(2,2)_T,(2,1)^2_T,(2,2)^2_T,(2,2)\right),$$
$$8\left(4,(2,2)^5_T,(2,2)\right),\ \ \ \ 15(7,(4,4)^5_T,(3,3)),$$
respectively.

To find $B'=D_4$ we proceed as follows. In an affine plane $\mathbb A,$ we fix a
smooth conic $C$ not containing the origin $p_0$ and we choose distinct points $p_1,\ldots,p_5\in C$
such that $T_1$ is tangent to $C$ at $p_1.$
We compute the linear system $L$ of plane curves of type
$15(7,(4,4)^5_T)$ and we resolve the base point of $L$ at $p_1,$
denoting the resulting linear system by $L_1.$
The defining polynomial of ${B'}$ is given by (the blow-down of)
a linear combination of elements of $L_1.$

In order to obtain a $(3,3)$-point $p_6$ (infinitely near to $p_1$),
we impose conditions, given by the zeros of
the derivatives up to order 2, to the elements of $L_1.$
Then we resolve this point and impose another triple
(infinitely near) point. With all these conditions we define a matrix, which is
denoted by $Mt$ in Appendix \ref{Appendix2}.
To have a solution it is necessary that $Mt$ has no maximal rank.

We define a scheme $Sch$ by imposing the following conditions:
\begin{description}
  \item $\cdot$ the vanishing of the maximal minors of $Mt$;
  \item $\cdot$ the points are infinitely near.
\end{description}
Now we compute the points of $Sch$ with Magma, choosing one of the solutions
for $p_6,$ and we use the procedures of Section \ref{Proc}
to compute ${B'},$ with singularities at $p_0,\ldots,p_6.$
Finally we perform some verifications: that ${B'}$ is reduced,
the singularities are as expected, etc.

Notice that there is no need to verify the non-existence of other singularities: 
in the presence of another non-negligible singularity, the computation of the invariants of the corresponding
double plane (using the formulas of \cite[Section V. 22]{BPV}) leads to a contradiction.\\

To find the pencil which induces the Albanese fibration,
one uses again the procedures $LinSys$ and $LinSys2$ to compute the pencil of curves of type
$$16(8,(4,4)_T^5,(4,4)),$$ through $p_0,\ldots,p_6,$ and to verify that the elements
$F_A^1,\ldots,F_A^4$ are as claimed.\\

With this we find a minimal double plane with $p_g=q=1,$ $K^2=7$ and $g=5.$
Using the above Magma procedures, one can verify that there
exist also branch loci of type
$$18\left(10,(5,5)^4_T,4,(3,3)\right)\ \ \ {\rm and}\ \ \ 16\left(8,(5,5)^3_T,4^2,(3,3)\right),$$
through $p_0,\ldots,p_6.$
These correspond to minimal Du Val double planes
with $p_g=q=1$ and $(K^2,g)=(5,4),(3,3).$
The pencils which induce the Albanese fibration are of type
$$15\left(7,(4,4)^4_T,4,(4,4)\right)\ \ \ {\rm and}\ \ \ 14\left(6,(4,4)^3_T,4^2,(4,4)\right),$$
respectively.\\

The Magma calculations for this section are in Appendix \ref{Appendix2}.
\subsection{$K^2=6$ and $g=3$}\label{cncTan}

In this section we impose two $(3,3)$-points to the branch locus of a Du Val's ancestor of type $\mathcal D_4,$
{\it i.e.} we construct a plane curve $B'$ of type $$14(6,(4,4)^4_T,(3,3)^2),$$
singular at points $p_0,\ldots,p_6,$ such that $B:=B'+\sum_1^4 T_i$ is reduced.
In this example the $(3,3)$-points $p_5,p_6$ are tangent to the conic $C$ through $p_1,\ldots,p_6.$
The double cover with branch locus $B$ is a plane model
of a Du Val double plane $S$ with $K^2=6,$ $p_g=q=1$ and $g=3.$

The pencil $f_A$ which induces the Albanese fibration of $S$ contains four fibres
$$F_A^1=D_1+2T_0,\ \ \ \ F_A^2=D_2,\ \ \ \ F_A^3=2D_3,\ \ \ \ F_A^4=2C+T_1+\cdots+T_4,$$
where $T_0$ is a line through $p_0,p_5,p_6$ and $D_1,D_2,D_3$ are curves of types
$$6\left(2,(2,2)^4_T,(1,1)^2\right),\ \ \ \ \ 8\left(4,(2,2)^4_T,(2,2)^2\right),\ \ \ \ 4\left(2,(1,1)^4_T,(1,1)^2\right),$$
respectively.
We have $$B'=D_1+D_2.$$
First we construct $D_1$ and then we use the procedure
$LinSys$ of Section \ref{Proc} to obtain $D_2$
(a general element of $f_A$) and to verify the existence of $D_3.$\\

To find $D_1,$ we follow the steps used in Section \ref{K8642}, but now with
$L$ the linear system of curves of type $6\left(2,(2,2)^4_T\right)$
and $Sch$ defined by
\begin{description}
  \item $\cdot$ $p_5,p_6\in C\bigcap F;$
  \item $\cdot$ $F$ is smooth at $p_5,p_6;$
  \item $\cdot$ $F$ is tangent to $C$ at $p_5,p_6;$
  \item $\cdot$ $p_5\ne p_6$ and $p_5,p_6\not\in\{p_0,\ldots,p_4\}.$
\end{description}

The details can be found in Appendix \ref{Appendix3}. Again we use
symmetry in order to increase speed of calculations.
\subsection{$K^2=8$ and $g=4$ or $3$}\label{K8g4}

Given, in an affine plane $\mathbb A,$ a general plane cubic $Q$ not containing the
origin $p_0,$ there are 6 points $p_1,\ldots,p_6\in Q$ such that each line $T_i,$
defined by $p_0,p_i,$ is tangent to $Q.$
In this example we have
$$F_A^1=3D_1,\ \ \ F_A^2=D_2,\ \ \ F_A^3=2D_3,\ \ \ F_A^4=2Q+T_1+\ldots+T_6,$$
where $D_1,D_2,D_3$ are curves of types
$$4\left(2,(1,1)^6_T\right),\ \ \ \ 12\left(6,(3,3)^6_T\right),\ \ \ \ 6\left(3,(2,1)^6_T\right),$$
respectively.

Using the procedure $LinSys$ (see Section \ref{Proc}), one can verify the
existence of the curves $D_1,D_2,D_3,$ through $p_0,\ldots,p_6,$ and that $p_1,\ldots,p_6$
are contained in a conic.
The double cover with branch locus $D_1+D_2+\sum_1^6T_i$ is a plane model
of a Du Val double plane $S$ with $K^2=6,$ $p_g=q=1$ and $g=4.$
This is a surface of type II described in \cite{Po4}.\\

Finally the surface of type I described in \cite{Po4}:\\
Let $C\subset\mathbb P^2$ be a smooth conic and $p_0\not\in C,$ $p_1,\ldots,p_6\in C$ be points such
that the lines $T_5$ and $T_6,$ through $p_0,p_5$ and $p_0,p_6,$ are tangent to $C.$
One can verify the existance of a smooth cubic $D$ of type $3(1,(1,1)_T^4,1^2),$
through $p_0,\ldots,p_6.$
Let $$F_A^1=2C+T_1+\cdots+T_4,\ \ \ \ \ \ F_A^2=2D+T_5+T_6$$ and $F_A^3,$ $F_A^4$
be general elements of the pencil generated by $F_A^1$ and $F_A^2.$ The branch locus in this case is
$$F_A^3+F_A^4+\sum_1^6T_i.$$
\subsection{$K^2=8,$ $g=4,$ $\phi_2$ not composed with $i$}\label{ExplBirrBic}

This section describes the steps to obtain an equation of a double plane $S$
of general type with $K^2=8,$ $p_g=q=1$ and Albanese fibres of genus $g=4,$
such that the bicanonical map 
$\phi_2$ of $S$ is not composed with the associated involution $i$ of $S.$

In this example the pencil $f_A$ which induces the Albanese fibration has four fibres
$$F_A^1=2D_1+T_5,\ \ \ F_A^2=6D_2+3T_1,\ \ \ F_A^3=2D_3+T_2+T_3+T_4,\ \ \ F_A^4=D_4,$$
where $D_1,\ldots,D_4$ are curves of types
$$10\left(4,(3,3)^4_T,(3,2)_T\right),\ \ \ \ \ \ 3\left(1,1,(1,1)^4_T\right),$$
$$9\left(3,(3,3),(3,2)^3_T,(3,3)\right),\ \ \ \ \ 21\left(9,(6,6)^5_T\right),$$
respectively. 
The surface $S$ is the minimal double plane given by the branch locus $$D_4+\sum_1^5T_i.$$
From the double cover formulas (\cite[Section V. 22]{BPV}), one obtains $K_S^2=8$
(notice that the pullback of the lines $T_i$ contains $10$ $(-1)$-curves), $\chi(\mathcal O_S)=1$
and $p_g(S)$ is equal to the number of generators of the linear system of curves of type $10(6,(2,3)_T^5).$

The points $p_0,\ldots,p_5$ are chosen such that the line $T_1,$ through $p_0,p_1,$ is tangent to the conic defined by $p_1,\ldots,p_5.$
One can verify that this implies the existence of the cubic $D_2.$ After finding $D_1,$ we verify the existence of $D_2$ and $D_3.$
The curve $D_4$ is a general element of the pencil generated by $D_1,D_2,D_3.$
One can use the procedure $LinSys$ to compute the number of elements of the linear system of curves of type
$10(6,(2,3)_T^5),$ showing that $p_g(S)=1.$

The difficulty here is the computation of $D_1.$ This is done following the steps of
Section \ref{K7}. We omit the details.
\appendix
\section{Appendix: Magma computations}\label{Magma}

In this appendix several computations are done using the {\em Computational
Algebra System MAGMA}.\\

The following functions will be useful:
\begin{verbatim}
function D(F,i);P:=Parent(F);
return Derivative(F,P.i);end function;

function D2(F,i,j);P:=Parent(F);
return Derivative(Derivative(F,P.i),P.j);end function;
\end{verbatim}

{\it i.e.} from now on
{\rm D(F,i)} means $\frac{\partial F}{\partial P.i}$ and
{\rm D2(F,i,j)} means $\frac{\partial^2 F}{\partial P.j \partial P.i}.$
\subsection{The procedures $LinSys$ and $LinSys2$}\label{Procedure}

\begin{verbatim}
procedure LinSys(A,d,p,m1,m2,td,~L)
  //p,m1,... are tuples of p_i,m1_i,...
  x:=A.1;y:=A.2;//The coordinates of A.
  L:=LinearSystem(LinearSystem(A,d),p,m1);
  for j:=1 to #m2 do
    if #Sections(L) eq 0 then break;end if;
    a:=p[j][1];b:=p[j][2];
    Bup:=[Evaluate(Sections(L)[i],y,(x-a)*y+b) div (x-a)^m1[j]:\
    i in [1..#Sections(L)]];//The strict transform of the
    //blown-up curves.
    L1:=LinearSystem(A,Bup);
    L2:=LinearSystem(L1,A![a,td[j]],m2[j]);//Imposing the
    //infinitely near singularity.
    if #Sections(L2) eq 0 then L:=L2;break;end if;
    Bdn:=[Evaluate((x-a)^m1[j]*Sections(L2)[i],y,(y-b)/(x-a)):\
    i in [1..#Sections(L2)]];//The blown-down curves.
    R:=Universe(Bdn);
    //R is a Rational function field. We need an homomorphism
    //to send the elements of Bdn into a polynomial ring.
    h:=hom<R->CoordinateRing(A)|[x,y]>;
    L:=LinearSystem(A,[h(Bdn[i]):i in [1..#Bdn]]);
  end for;
end procedure;
\end{verbatim}

\begin{verbatim}
procedure LinSys2(A,L,q,m,td,~J)
  x:=A.1;y:=A.2;//The coordinates of A.
  J:=LinearSystem(L,q,m[1]);
  td:=[q[2]] cat td;
  for j:=1 to #td-1 do
    if #Sections(J) eq 0 then break;end if;
    b:=td[j];
    Bup:=[Evaluate(Sections(J)[i],y,(x-q[1])*y+b) div\
    (x-q[1])^m[j]:i in [1..#Sections(J)]];
    J1:=LinearSystem(A,Bup);
    J:=LinearSystem(J1,A![q[1],td[j+1]],m[j+1]);
  end for;
  //
  for j:=#td-1 to 1 by -1 do
    if #Sections(J) eq 0 then break;end if;
    b:=td[j];
    Bdn:=[Evaluate((x-q[1])^m[j]*Sections(J)[i],y,(y-b)/\
    (x-q[1])):i in [1..#Sections(J)]];
    R:=Universe(Bdn);
    h:=hom<R->CoordinateRing(A)|[x,y]>;
    J:=LinearSystem(A,[h(Bdn[i]):i in [1..#Bdn]]);
  end for;
end procedure;
\end{verbatim}
\subsection{Double planes}\label{dblplns}

\subsubsection{$K^2=8,6,4,2$ and $g=5,4,3,2$}\label{Appendix1}
Here we have the computations of Section \ref{K8642}, with details
for the case of a double plane with $p_g=q=1$ and $K^2=6.$
\begin{verbatim}
> A<x,y>:=AffineSpace(Rationals(),2);
> p:=[A![7/5,4/5],A![7/5,-4/5],A![2,1],A![2,-1],Origin(A)];
> d:=6;m1:=[2,2,2,2,2];m2:=[2,2,1,1];
> td:=[p[i][2]/p[i][1]:i in [1..#m2]];
> LinSys(A,d,p,m1,m2,td,~L);
> #Sections(L);BaseComponent(L);
5
Scheme over Rational Field defined by
1
\end{verbatim}
We are using symmetry: we want to find $p_5=(x,y)$ and $p_6=(x,-y).$
\begin{verbatim}
> R<x,y,b,c,d,e,n>:=PolynomialRing(Rationals(),7);
> h:=hom<PolynomialRing(L)->R|[x,y]>;
> l:=h(Sections(L));
> F:=l[1]+b*l[2]+c*l[3]+d*l[4]+e*l[5];
> G:=Evaluate(F,y,-y);
> eqF:=2*x*y*D2(F,1,2)+x^2*D2(F,1,1)+y^2*D2(F,2,2);
> eqG:=Evaluate(eqF,y,-y);
\end{verbatim}
The condition $eqF=0$ is used to obtain one branch of the
double point $p_5=(x,y)$ tangent to the line $T_5.$
In fact, the tangent cone of a plane curve $\{F(X,Y)=0\}$ at $(x,y)$ is given by
$$\frac{\partial^2 F}{\partial X^2}(x,y)\frac{(X-x)^2}{2}+\frac{\partial^2 F}{\partial X\partial Y}(x,y)(X-x)(Y-y)+
\frac{\partial^2 F}{\partial Y^2}(x,y)\frac{(Y-y)^2}{2}.$$
One can verify that this equation is divisible by $xY-yX$ only if 
$$2xy\frac{\partial^2 F}{\partial X\partial Y}(x,y)+x^2\frac{\partial^2 F}{\partial X^2}(x,y)+
y^2\frac{\partial^2 F}{\partial Y^2}(x,y)=0.$$
\begin{verbatim}
> dif:=y*(y-1/2*x)*(y-4/7*x)*(y+1/2*x)*(y+4/7*x);
\end{verbatim}
In order to obtain $p_5,p_6\not\in T_i,$ $i=1,\ldots,4,$ and $p_5\not=p_6,$
we need $dif$ to be {\em different} from zero.
This is achieved by imposing the condition $1+n\cdot dif=0.$
\begin{verbatim}
> A:=AffineSpace(R);
> Sch:=Scheme(A,[(x-2)^2+y^2-1,F,D(F,1),D(F,2),G,D(G,1),D(G,2),\
> eqF,eqG,1+n*dif]);
> Dimension(Sch);
0
> PointsOverSplittingField(Sch);
\end{verbatim}

This command gives the points of $Sch$ and the necessary field extensions
to define them. Choosing one of the solutions we obtain points $p_5,$ $p_6$
such that there exists a pencil $f_A$ of curves of type $15(7,(4,4)^5_T,4).$
Let $B'$ be a general element of $f_A.$
The branch locus $B:=B'+\sum_1^5 T_i$ is of type $20(12,(5,5)^5_T,4).$
The corresponding minimal double plane is a surface of general type with $p_g=q=1,$
$K^2=6$ and $g=4.$
\\

Verification that $B_0:=B'$ is as stated:

\begin{verbatim}
> R<r3>:=PolynomialRing(Rationals());
> K<r3>:=NumberField(r3^4 - 570063504574501/8986626*r3^2+\
> 194676993199491455085153141001/323037787455504);
> x1:=-1225449/218906496039245*r3^2 + 6763320857703/401161\
> 4254780;
> y1:=-7879209182423568/1971150953143623770162761495*r3^3 +
> 37270947258282841632/117281546566527266624785*r3;
> //
> A<x,y>:=AffineSpace(K,2);
> p:=[A![7/5,4/5],A![7/5,-4/5],A![2,1],A![2,-1],A![x1,y1],\
> A![x1,-y1],Origin(A)];
> d:=15;m1:=[4,4,4,4,4,4,7];m2:=[4,4,4,4,4];
> td:=[p[i][2]/p[i][1]:i in [1..#m2]];
> LinSys(A,d,p,m1,m2,td,~L);
> #Sections(L);BaseComponent(L);
2
Scheme over K defined by
1
> B0:=Curve(A,Sections(L)[1]+Sections(L)[2]);
> IsReduced(B0);
true
> Multiplicity(B0,Origin(A));
7
> IsOrdinarySingularity(B0,Origin(A));
true
> [Multiplicity(B0,p[i]):i in [1..6]];
[ 4, 4, 4, 4, 4, 4 ]
> IsOrdinarySingularity(B0,p[6]);
true
> T:=[Curve(A,y-p[i][2]/p[i][1]*x):i in [1..5]];
> [IntersectionNumber(B0,T[i],p[i]):i in [1..5]];
[ 8, 8, 8, 8, 8 ]
> ResolutionGraph(B0,p[1]);
The resolution graph on the Digraph
Vertex      Neighbours

1 ([ -2, 4, 1, 0 ]) 2 ;
2 ([ -1, 8, 2, 4 ]) ;
\end{verbatim}
One obtains the same resolution graph for $p_2,...,p_5.$\\

The singularities of a general element of $L$ are no worst than the ones above,
hence we do not need to verify that $B_0$ has no other singularities.\\

The other cases, $p_g=q=1,$ $K^2=8,4,2$ and $g=5,3,2,$ are analogous to the previous one.
One needs only to ask Magma (using the procedure $LinSys$) for curves of type
$$[10+2n](2n+2,(5,5)^n_T,4^{6-n}),\ \ \ n=6,4,3,$$
with singularities at the previous points $p_0,\ldots,p_6.$
\subsubsection{$K^2=7,5,3$ and $g=5,4,3$}\label{Appendix2}
The detailed computations of Section \ref{K7} are as follows.
\begin{verbatim}
> A<x,y>:=AffineSpace(Rationals(),2);
> p:=[A![1,0],A![1/5,2/5],A![2/5,1/5],A![8/5,9/5],\
> A![9/5,8/5],Origin(A)];
> d:=15;m1:=[4,4,4,4,4,7];m2:=[4,4,4,4,4];
> td:=[p[i][2]/p[i][1]:i in [1..#m2]];
> LinSys(A,d,p,m1,m2,td,~L);
> #Sections(L);BaseComponent(L);
8
Scheme over Rational Field defined by
1
> Bup:=[Evaluate(Sections(L)[i],y,(x-1)*y+0) div (x-1)^4:\
> i in [1..#Sections(L)]];
> L:=LinearSystem(A,Bup);
> Bup:=[Evaluate(Sections(L)[i],y,(x-1)*y+0) div (x-1)^4:\
> i in [1..#Sections(L)]];
> L1:=LinearSystem(A,Bup);
\end{verbatim}
At this stage we have imposed the $(4,4)$-points and resolved $p_1.$
\begin{verbatim}
> R<x,y,u,v,n>:=PolynomialRing(Rationals(),5);
> h:=hom<PolynomialRing(L1)->R|[x,y]>;
> l:=h(Sections(L1));
\end{verbatim}
Now we impose, to the elements of $l,$ the necessary conditions
in order to obtain the $(3,3)$-point $p_6=(u,v).$
The matrix $Mt$ defined by these conditions cannot have maximal rank.
\begin{verbatim}
> H:=[Evaluate(l[i],[u,v,u,v,n]):i in [1..#l]];
> F:=[Evaluate(l[i],y,(x-u)*y+v):i in [1..#l]];
> G:=[(F[i]-Evaluate(F[i],x,u)) div (x-u):i in [1..#l]];
> G1:=[(G[i]-Evaluate(G[i],x,u)) div (x-u):i in [1..#l]];
> G2:=[(G1[i]-Evaluate(G1[i],x,u)) div (x-u):i in [1..#l]];
> F:=G2;
> M:=[[H[i],D(H[i],3),D(H[i],4),D2(H[i],3,3),D2(H[i],3,4),\
> D2(H[i],4,4),F[i],D(F[i],1),D(F[i],2),D2(F[i],1,1),\
> D2(F[i],1,2),D2(F[i],2,2)]:i in [1..#l]];
> ME:=[[Evaluate(M[i][o],[1,y,1,v,n]):o in [1..12]]:i in [1..#M]];
> //This last step is needed to increase the speed of calculations.
> Mt:=Matrix(ME);
> min:=Minors(Mt,#l);
> A:=AffineSpace(R);
> Sch:=Scheme(A,min cat [x-1,u-1,1+n*v]);
> Dimension(Sch);
0
> PointsOverSplittingField(Sch);

\end{verbatim}
As before, this gives various solutions. We choose one who works. Here goes
the verifications:
\begin{verbatim}
> R<r2>:=PolynomialRing(Rationals());
> K<r2>:=NumberField(r2^2 - 1292/35);
> A<x,y>:=AffineSpace(K,2);
> p:=[A![1/5,2/5],A![2/5,1/5],A![8/5,9/5],A![9/5,8/5],\
> Origin(A)];
> d:=15;m1:=[4,4,4,4,7];m2:=[4,4,4,4];
> td:=[p[i][2]/p[i][1]:i in [1..#m2]];
> LinSys(A,d,p,m1,m2,td,~L);
> //
> q:=A![1,0];m:=[4,4,3,3];
> td:=[0,35/1292*r2,-455/15504*r2 + 455/7752];
> LinSys2(A,L,q,m,td,~J);
> #Sections(J);
1
> B0:=Curve(A,Sections(J)[1]);
> IsReduced(B0);
true
> p1:=[q] cat p;
> T:=[Curve(A,y-p1[i][2]/p1[i][1]*x):i in [1..5]];
> [Multiplicity(B0,p1[i]):i in [1..6]];
[ 4, 4, 4, 4, 4, 7 ]
> [IntersectionNumber(T[i],B0,p1[i]):i in [1..5]];
[ 8, 8, 8, 8, 8 ]
> [ResolutionGraph(B0,p1[i]):i in [1,2]];
[
    The resolution graph on the Digraph
    Vertex      Neighbours

    1 ([ -2, 8, 2, 1 ]) 2 3 ;
    2 ([ -2, 4, 1, 0 ]) ;
    3 ([ -2, 11, 3, 0 ])        4 ;
    4 ([ -1, 14, 4, 3 ])        ;
    ,
    The resolution graph on the Digraph
    Vertex      Neighbours

    1 ([ -2, 4, 1, 0 ]) 2 ;
    2 ([ -1, 8, 2, 4 ]) ;
]
\end{verbatim}
The resolution graphs for the points $p_3,$ $p_4$ and $p_5$ are equal
to this last one. \\

Now we calculate the pencil $J$ which induces the Albanese fibration.
\begin{verbatim}
> d:=16;m1:=[4,4,4,4,8];m2:=[4,4,4,4];
> td:=[p[i][2]/p[i][1]:i in [1..#m2]];
> LinSys(A,d,p,m1,m2,td,~L);
> q:=A![1,0];m:=[4,4,4,4];
> td:=[0,35/1292*r2,-455/15504*r2 + 455/7752];
> LinSys2(A,L,q,m,td,~J);
> #Sections(J);BaseComponent(J);
2
Scheme over K defined by
1
> Jy:=LinearSystem(J,Curve(A,y));
> B0 eq Curve(A,Sections(Jy)[1] div y);
true
\end{verbatim}
With this we have constructed a minimal double plane with $p_g=q=1,$ $K^2=7$ and $g=5.$
\subsubsection{$K^2=6$ and $g=3$}\label{Appendix3}
Here we give the detailed computations of Section \ref{cncTan}.
\begin{verbatim}
> A<x,y>:=AffineSpace(Rationals(),2);
> p:=[A![4/5,7/5],A![-4/5,7/5],A![1,2],A![-1,2],Origin(A)];
> d:=6;m1:=[2,2,2,2,2];m2:=[2,2,2,2];
> td:=[p[i][2]/p[i][1]:i in [1..#m2]];
> LinSys(A,d,p,m1,m2,td,~L);
> #Sections(L);BaseComponent(L);
2
Scheme over Rational Field defined by
1
> R<x,y,b,n>:=PolynomialRing(Rationals(),4);
> h:=hom<PolynomialRing(L)->R|[x,y]>;
> l:=h(Sections(L));
> F:=l[1]+b*l[2];
> G:=Evaluate(F,x,-x);
> C:=x^2+(y-2)^2-1;
> //
> eqF:=D(C,1)*D(F,2)-D(C,2)*D(F,1);//To obtain a curve
> //tangent to the conic C at p_5=(x,y).
> eqG:=Evaluate(eqF,x,-x);//The same to p_6=(-x,y).
> dif:=x*(y-2*x)*(y-7/4*x)*(y+2*x)*(y+7/4*x);
> //We need dif to be non-zero.
> //
> A:=AffineSpace(R);
> Sch:=Scheme(A,[C,F,G,eqF,eqG,1+n*D(F,1)*D(F,2)*dif]);
> Dimension(Sch);
0
> PointsOverSplittingField(Sch);
{@(-r2, 0, 5377/5292, -3/307328), (r2, 0, 5377/5292, -3/307328)@}
Algebraically closed field with 2 variables
Defining relations: [
    r2^2 + 3,
    r1^2 + 3 ]
\end{verbatim}

This gives the points $p_5,$ $p_6.$
We omit the remaining verifications (they are similar to the ones in \ref{Appendix1} and \ref{Appendix2}).
We notice only that, to verify that a curve $D$ has no singularity at infinity, one can proceed as follows:

\begin{verbatim}
> PD:=ProjectiveClosure(D);
> Dimension(SingularSubscheme(PD) meet LineAtInfinity(A));
\end{verbatim}
\bibliography{ReferencesRito1}

\bigskip
\bigskip

\noindent Carlos Rito
\\ Departamento de Matem\' atica
\\ Universidade de Tr\' as-os-Montes e Alto Douro
\\ 5000-911 Vila Real
\\ Portugal
\\\\
\noindent {\it e-mail:} crito@utad.pt

\end{document}